\documentclass[11pt,a4paper,reqno]{amsart}
\usepackage{graphicx,amsmath,latexsym,amssymb,amsthm,mathrsfs, amscd, amsfonts, color}
\usepackage{filecontents}
\def\MR#1{}

\usepackage{enumerate}
\def\MR#1{}
\newtheorem{theorem}{Theorem}[section]
\newtheorem{lemma}[theorem]{Lemma}

\theoremstyle{definition}
\newtheorem{definition}[theorem]{Definition}
\newtheorem{example}[theorem]{Example}

\numberwithin{equation}{section}


\begin{document}
\noindent {\footnotesize\tiny}\\[1.00in]
\textcolor[rgb]{0.00,0.00,1.00}{}

\title[A Mean Ergodic Theorem in Hadamard spaces]{A Mean Ergodic Theorem for Nonexpansive Mappings in Hadamard spaces}
\maketitle
\begin{center}
	{\sf Hadi Khatibzadeh$^{1}$ and Hadi Pouladi$^2$}\\
	{\footnotesize{\it $^{1,2}$ Department of Mathematics, University
			of Zanjan, P. O. Box 45195-313, Zanjan, Iran.\\ Email:$^1$hkhatibzadeh@znu.ac.ir, $^{2}$hadi.pouladi@znu.ac.ir.}}
\end{center}

\begin{abstract}
In this paper, we prove a mean ergodic theorem for nonexpansive mappings in Hadamard (nonpositive curvature metric) spaces, which extends the Baillon nonlinear ergodic theorem. The main result shows that the sequence given by the Karcher means of iterations of a nonexpansive mapping with a nonempty fixed point set converges weakly to a fixed point of the mapping. This result also remains true for a $1-$parameter continuous semigroup of contractions.
\end{abstract}
\section{Introduction and Preliminaries}
Let $H$ be a real Hilbert space. A mapping $T:H\longrightarrow H$ is called nonexpansive if\linebreak $\|Tx -Ty\|\leqslant \|x -y\|$, for all $x,y\in H$. $F(T)=\{x : Tx=x\}$ denotes the set of all fixed points of the mapping $T$. The sequence $\{T^nx\}$ of iterations of a mapping $T$ at a point $x$ is called strongly (weakly) Cesaro convergent if $\frac{1}{n}\sum_{i=0}^{n-1} T^ix$ converges strongly (weakly) to a point. The study of mean ergodic theorems began with von Neumann in 1932 \cite{VNeumann1932} for linear nonexpansive mappings in Hilbert spaces. Birkhoff in 1939 extended this theorem to Banach spaces \cite{Birkhoff1939}. About forty years later in 1975 Baillon \cite{baillon1975} considered the nonlinear version of the mean ergodic theorem of von Neumann that is well-known as Baillon's  nonlinear mean ergodic theorem in the literature as follow:
\begin{theorem}\label{baillon}
Let $H$ be a Hilbert space and $C$ be a nonempty, closed and convex subset of $H$. Let $T:C\longrightarrow C$ be a nonexpansive mapping with $F(T)\neq \emptyset$. Then, for any $x\in C$ the Cesaro mean $\frac{1}{n}\sum_{k=0}^{n-1}T^kx$ converges weakly to a fixed point of $T$ as $n\rightarrow\infty$.
\end{theorem}
\noindent Br\'{e}zis and Browder \cite{Brezis1976} extended Baillon's nonlinear ergodic theorem to convergence of more general summability methods. Reich \cite{Reich1978} improved these results and simplified their proof and in \cite{Reich1979} generalized these results to Banach spaces. Bruck \cite{Bruck1979} provided another proof for the ergodic theorem in Banach spaces. After them, many authors and researchers studied various versions of nonlinear ergodic theorems for nonexpansive mappings and semigroups and their generalizations in linear spaces setting mainly in Banach spaces. The reader can consult \cite{takahashi2007} and \cite{takahashi2014} and references there in to see a relatively complete bibliography and a list of works in this field.

Kirk \cite{Kirk2003} and \cite{Kirk2004} with studying of fixed point theory in the setting of Hadamard spaces showed that they are a suitable framework for studying nonlinear analysis. After him, extensive studies were conducted to nonexpansive mappings and fixed point theory in Hadamard spaces. In \cite{kakavandi2015} and \cite{kakavandi2011} Ahmadi Kakavandi studied Ballion's mean ergodic theorem for a nonlinear amenable semigroup of nonexpansive mappings. In \cite{kakavandi2015} he showed ergodic convergence of the semigroup but he assumed that the mean satisfies a certain property, namely,  $\mathcal{L}$ {\it property} (see \cite{kakavandi2015}), which is not satisfied even by usual means. Liimatainen \cite{Liimatainen2012} proved a result related to the mean ergodic convergence for orbits of nonexpansive mappings. He showed the strong convergence of orbits of nonexpansive mappings which have an additional condition, namely, {\it distance convexity}. This condition is an  extension of linearity in Hadamard spaces and his result can be considered as an extension of von Neumann mean ergodic theorem in this setting. But for general nonexpansive mappings (without additional assumptions), the main result of \cite{Liimatainen2012} implies that every weak cluster point of the mean of a bounded  orbit is a fixed point of the mapping. The mean ergodic theorem for the orbit of a general nonexpansive mapping is still an open problem in general Hadamard spaces, because the set of weak cluster points is not necessarily a singleton. In this paper we establish the mean ergodic theorem for general nonexpansive mappings in Hadamard spaces with an additional condition, namely, $(\overline{Q_4})$ condition (defined blow in this section), as:
\begin{theorem}\label{weakergodickarcher}
Let $(\mathscr{H},d)$ be a Hadamard space with $(\overline{Q_4})$ condition, $C$ be a nonempty, closed and convex subset of $\mathscr{H}$ and $T:C\rightarrow C$ be a nonexpansive mapping with $F(T)\neq\emptyset$. Then the sequence $\{T^nx\}$ is $\triangle -mean~convergent$ to a fixed point of $T$, which is also the strong limit point of the sequence $\{PT^nx\}$.
\end{theorem}
First, we briefly present some definitions and preliminaries. Let $(X,d)$ be a metric space. A geodesic segment (or geodesic) between two points $x_0,x_1\in X$, is the image of an isometry mapping $\gamma:[0, d(x_0,x_1)]\longrightarrow X$, with $\gamma(0)=x_0, \gamma(d(x_0,x_1))=x_1$ and $d\big(\gamma(t),\gamma(t')\big)=|t-t'|~\text{for all}~t,t'\in[0,d(x,y)]$.
A metric space $(X,d)$ is said to be a geodesic metric space if every two points of $X$ are joined by a geodesic and it said to be uniquely geodesic if between any two points there is exactly one geodesic that for two arbitrary points $x_0,x_1$ is denoted by $[x_0,x_1]$. All points in $[x_0,x_1]$ are denoted by $x_t=(1-t)x_0\oplus tx_1$ for all $t\in[0,1]$, where $d(x_t,x_0)=td(x_0,x_1)$ and $d(x_t,x_1)=(1-t)d(x_0,x_1)$.
In the uniquely geodesic space $X$, $C\subset X$ is said to be convex if for each $x,y\in C$, $[x,y]\subset C$. For an arbitrary subset $A$ of $X$, the closed convex hull of $A$ is the smallest closed convex set that contians $A$. We denote the closed convex hull of $A$ by $\overline{co}(A)$.
\noindent A function $f: X\longrightarrow \Bbb{R}$ on a uniquely geodesic metric space $(X,d)$ is said to be convex if for all $x,y\in X$ and all $\lambda\in [0,1]$,
\begin{equation*}
f\big((1-\lambda) x\oplus\lambda y\big)\leqslant (1-\lambda)f(x)+\lambda f(y).
\end{equation*}
$f$ is said to be strongly convex with parameter $\gamma>0$ if for all $x,y\in X$ and all $\lambda\in [0,1]$,
\begin{equation*}
f\big((1-\lambda) x\oplus \lambda y\big)\leqslant (1-\lambda) f(x)+\lambda f(y)-\lambda(1-\lambda)\gamma d^2(x,y).
\end{equation*}
As proved in \cite[Lemma 2.5]{DHOMPONGSA20082572} and \cite[page 163]{bridson2011metric}, a uniquely geodesic metric space $(X,d)$ is a $CAT(0)$ space if and only if, for every $x\in X$, the function $d^2(x,\cdot)$ is strongly convex with $\gamma=1$. In other words, for every three points $x_0 ,x_1, y \in X$ and for every $0<t<1$,
\allowdisplaybreaks\begin{equation*}
d^2(y, x_t)\leqslant(1-t)d^2(y,x_0)+td^2(y,x_1)-t(1-t)d^2(x_0,x_1)
\end{equation*}
is satisfied, where $x_t=(1-t)x_0\oplus tx_1$ for every $t\in[0,1]$.

\noindent A complete $CAT(0)$ space is said to be a Hadamard space. From now, we denote every Hadamard space by $\mathscr H$. A function $f:\mathscr H\longrightarrow \mathbb{R}$ is said to be lower semicontinuous if the set $\{x\in\mathscr H : f(x)\leqslant \alpha\}$ is closed for all $\alpha\in \Bbb{R}$. Any lower semicontinuous, strongly convex function on a Hadamard space has a unique minimizer \cite[Proposition 2.2.17]{bacak2014convex}.

For $a,b,c,d\in X$, we denote
\begin{equation*}\label{e12}
\frac{1}{2}\big\{d^2(a,d)+d^2(b,c)-d^2(a,c)-d^2(b,d) \big\}
\end{equation*}
by $\langle ab,cd \rangle$, which is called quasi-inner product. Berg and Nikolaev in \cite[Corollary 3]{Berg2008} proved that $CAT(0)$ spaces satisfy the Cauchy-Schwarz like inequality:
\begin{equation*}
\langle ab, cd\rangle\leqslant d(a,b)d(c,d), \quad \forall a,b,c,d\in X.
\end{equation*}
\noindent The following lemma is easily proved.
\begin{lemma}\label{l8}
Let $(X,d)$ be a $CAT(0)$ space and $a,b,c,d,e\in X$. Then,
\begin{enumerate}[(i)]
\item $\langle ab,cd \rangle=\langle cd,ab \rangle$.
\item $\langle ab,cd \rangle=-\langle ab,dc \rangle=-\langle ba,cd \rangle$.
\item $\langle ab,cd \rangle=\langle ae,cd \rangle+\langle eb,cd \rangle$.
\end{enumerate}
\end{lemma}
A mapping $T:X\longrightarrow X$ on a metric space $(X,d)$ is called nonexpansive if\linebreak $d(Tx,Ty)\leqslant d(x,y)$, for all $x,y\in X$. As proved in \cite[Theorem 12]{Kirk2004}, for a mapping $T$ on a Hadamard space, $F(T)$ is closed and convex. Also, in Hadamard spaces for any nonempty closed convex subset $S$, $P_Sx:=\{s\in S: d(x,S)=d(x,s) \}$ is a singleton, where $d(x,S):=\inf_{s\in S} d(x,s)$\cite[Theorem 2.1.12]{bacak2014convex}. Thus, the metric projection on nonempty closed convex subset $S$ is the following map:
\begin{equation*}
\underset{x\longmapsto P_Sx}{P:\mathscr H\longrightarrow S},
\end{equation*}
where $P_Sx$ is the nearest point of $S$ to $x$ for all $x\in \mathscr H$. A well-known fact implies that
\begin{equation}\label{projection}
d^2(x,P_Sx)+d^2(P_Sx,y)\leqslant d^2(x,y),\ \forall y\in S
\end{equation}
(see \cite[Theorem 2.1.12]{bacak2014convex}).
For more facts about Hadamard spaces, the readers can consult \cite{bacak2014convex} and \cite{bridson2011metric}. 

\noindent Let $(\mathscr H,d)$ be a Hadamard space, $\{x_n\}$ be a sequence in $\mathscr H$ and $x\in \mathscr H$. Set
\begin{equation*}
r(x,\{x_n\})=\displaystyle\limsup_{n\rightarrow\infty} d(x,x_n).
\end{equation*}
The asymptotic radius of $\{x_n\}$ is defined as
\begin{equation*}
r(\{x_n\})=\displaystyle\inf_{x\in\mathscr H} r(x,\{x_n\}),
\end{equation*}
and the asymptotic center of $\{x_n\}$ is the set
\begin{equation*}
A(\{x_n\})=\big\{x\in \mathscr H : r(\{x_n\})=r(x,\{x_n\}\big\}.
\end{equation*}
A well-known result implies that  $A(\{x_n\})$ is a singleton in Hadamard spaces \cite[Proposition 7]{Dhompongsa2006}. 

\begin{definition}\label{weakcondef1}
A sequence $\{x_n\}$ in a Hadamard space $\mathscr H$ $\triangle -converges$ to $x\in \mathscr H$ if $A(\{x_{n_k}\})=\{x\}$ for every subsequence $\{x_{n_k}\}$ of $\{x_n\}$. In this case we write $\bigtriangleup-\lim_n x_n=x$ or $x_n\overset{\bigtriangleup}{\longrightarrow}x$.
\end{definition}
The notion of $\triangle-$convergence first introduced by Lim \cite{lim1976} in general metric spaces. Also, Kirk and Panyanak in \cite{KIRK20083689} showed that in Hadamard spaces this concept of convergence shares many properties of the weak convergence in linear spaces. 

\begin{lemma}$($see \cite{KIRK20083689}$)$.\label{lconsubseq}
Every bounded sequence in a Hadamard space has a $\triangle -convergent$ subsequence. Also every closed convex subset of a Hadamard space is $\bigtriangleup-$closed in the sense that it contains all $\bigtriangleup-limit$ point of every $\bigtriangleup-convergent$ sequence of points of the subset.
\end{lemma}

The following geometric condition for nonpositive curvature metric spaces has been introduced by Kirk and Panyanak \cite[page 3693]{KIRK20083689}:

$(\mathbf{Q_4})$ for points $x,y,p,q\in\mathscr H$ and any point $m$ in the segment $[x,y]$,
\begin{equation*}
d(p,x)<d(x,q)\ \&\ d(p,y)<d(y,q) \Longrightarrow d(p,m)\leqslant d(m,q).
\end{equation*}
Also the following modification of $(\mathbf{Q_4})$ condition was introduced by Kakavandi \cite[page 6]{kakavandi2013}:

$(\mathbf{\overline{Q_4}})$ for points $x,y,p,q\in\mathscr H$ and any point $m$ in the segment $[x,y]$,
\begin{equation*}
d(p,x)\leqslant d(x,q)\ \&\ d(p,y)\leqslant d(y,q) \Longrightarrow d(p,m)\leqslant d(m,q).
\end{equation*}
Hilbert spaces, $\Bbb{R}-$trees and any $CAT(0)$ space of constant curvature satisfy $(\mathbf{\overline{Q_4}})$ condition. Clearly $(\mathbf{\overline{Q_4}})$ implies $(\mathbf{Q_4})$ (see Kakavandi \cite{kakavandi2013}). Also $(\mathbf{\overline{Q_4}})$ condition implies that $F(x,y):=\{z\in \mathscr H\ : \ d(x,z)\leqslant d(z,y) \}$ is convex for any $x,y\in \mathscr H$.

\begin{lemma}[Demiclosedness in $CAT(0)$ space]$($see \cite{KIRK20083689}$)$.\label{demiclosedness}
Let $C$ be a closed convex subset of a $CAT(0)$ space $X$ and $T:C\longrightarrow X$ be a nonexpansive mapping. If $x_n$ is $\triangle -convergent$ to $x$ and $d(x_n,Tx_n)\longrightarrow 0$ , then $x\in F(T)$.
\end{lemma}

One of means defined in Hadamard space is the Karcher mean that is extension of the linear mean in Hilbert spaces. We state it for a sequence as follows.
\begin{definition}[Karcher mean]\label{Karcher mean}
Given a sequence $\{x_n\}_{n=0}^{+\infty}$  in a Hadamard space. For the $n$ first terms $x_0,\ldots ,x_{n-1}$ of the sequence, we define 
\begin{equation}\label{e22}
{\mathcal{F}}_n(x)=\frac{1}{n}\displaystyle\sum_{i=0}^{n-1} d^2(x_i,x),
\end{equation}
and for $x_k, \cdots, x_{k+n-1}$, $k\geqslant 1$, we define
\begin{equation}\label{e23}
\mathcal{F}^k_n(x)=\frac{1}{n}\displaystyle\sum_{i=0}^{n-1} d^2(x_{k+i},x).
\end{equation}
By \cite[Proposition 2.2.17]{bacak2014convex} these functions are lower semicontinuous and strongly convex and hence they have unique minimizers. For ${\mathcal{F}}_n(x)$ the unique minimizer is denoted by $\sigma_n(x_0,\ldots ,x_{n-1})$ (or briefly $\sigma_n$) and it is called the Karcher mean of $x_0,\ldots ,x_{n-1}$ that introduced in \cite{karcher1977}.  $\{\sigma_n\}$ is called the sequence of means of $\{x_n\}$. Also for the function $\mathcal{F}_n^k(x)$ the unique minimizer is denoted by $\sigma_n^k (x_k,\ldots ,x_{k+n-1})$ (or briefly $\sigma_n^k$), which is the Karcher mean of $x_k,\ldots ,x_{k+n-1}$.

\noindent Let $T:\mathscr H\to \mathscr H$ be a mapping. Fixed $p\in \mathscr H$. For the orbit $\{T^n p | n=0,1,2,\ldots\}$, $\sigma_n(p)$ and $\sigma_n^k(p)$, are defined respectively as the unique minimizers of the functions
\begin{equation*}
{\mathcal{F}}[p]_n(x)=\frac{1}{n}\displaystyle\sum_{i=0}^{n-1} d^2(T^ip,x),
\end{equation*}
and
\begin{equation*}
\mathcal{F}[p]_n^k(x)=\frac{1}{n}\displaystyle\sum_{i=0}^{n-1} d^2(T^{k+i}p,x),
\end{equation*}
i.e., $\sigma_n(p)$ is the Karcher mean of $p, Tp,\ldots, T^{n-1}p$ and $\sigma_n^k(p)$ is the Karcher mean of $T^kp, T^{k+1}p,\ldots, T^{k+n-1}p$.
\end{definition}
\section{Weak Ergodic Theorem}
In this section we prove the weak convergence of the sequence of the Karcher means $\{\sigma_n(x)\}$ defined in Definition \ref{Karcher mean} for the orbit of a nonlinear nonexpansive mapping in Hadamard spaces with $(\mathbf{\overline{Q_4}})$ condition. This result generalizes Baillon's mean ergodic theorem to a more general setting. We first recall the following lemma, which is  a consequence of \cite[Proposition 4.1]{kakavandi2011}.


\begin{lemma}\label{the2}
Let $\mathscr H$ be a Hadamard space and $T:\mathscr H\longrightarrow \mathscr H$ a nonexpansive mapping such that $F(T)$ is nonempty. Let $P$ be the metric projection from $\mathscr H$ onto $F(T)$. Then for any $x\in \mathscr H$, $\{PT^n x\}$ converges strongly to an element $p$ of $F(T)$. Moreover, $p$ is the unique asymptotic center of the orbit $\{T^nx\}$.
\end{lemma}
Following lemmas are needed to prove the main result.
\begin{lemma}\label{l15}
Let $\{x_n\}$ be a sequence in Hadamard space $\mathscr H$. Then for $\sigma_n^k$ defined as the above, for each $y\in\mathscr H$ and $k\geqslant 1$ we have:
\begin{enumerate}[(i)]
\item\label{l15i1}
$\displaystyle d^2\big( \sigma_n^k ,y\big)\leqslant \frac{1}{n}\displaystyle\sum_{i=0}^{n-1} d^2(x_{k+i},y)-\frac{1}{n}\displaystyle\sum_{i=0}^{n-1} d^2(x_{k+i},\sigma_n^k).$
\item\label{l15i2} $\displaystyle d\big( \sigma_n^k ,y\big)\leqslant \frac{1}{n}\displaystyle\sum_{i=0}^{n-1} d(x_{k+i},y).$
\end{enumerate}
\end{lemma}
\begin{proof}
\eqref{l15i1}. This part is a consequence of \cite[Lemma 2.7]{kakavandi2015} that we summarize its proof in this setting. Since $\sigma_n^k$ is the unique minimizer of $\mathcal{F}_n^k(x)$ defined in \eqref{e23} and by the strong convexity of this function, for $0<\lambda <1$ we have:
\allowdisplaybreaks\begin{eqnarray}
\mathcal{F}_n^k(\sigma_n^k)&\leqslant & \mathcal{F}_n^k\big(\lambda\sigma_n^k\oplus (1-\lambda)y\big) \nonumber \\
&\leqslant &\lambda\mathcal{F}_n^k(\sigma_n^k)+ (1-\lambda)\mathcal{F}_n^k(y)- \lambda(1-\lambda)d^2\big( \sigma_n^k ,y\big). \nonumber
\end{eqnarray}
Therefor we obtain
\begin{equation*}
\lambda d^2\big( \sigma_n^k ,y\big)\leqslant \mathcal{F}_n^k(y)-\mathcal{F}_n^k(\sigma_n^k).
\end{equation*}
Letting $\lambda\rightarrow 1$ implies:
\begin{eqnarray}\label{ee1}
d^2\big( \sigma_n^k ,y\big)&\leqslant &\mathcal{F}_n^k(y)-\mathcal{F}_n^k(\sigma_n^k) \nonumber \\
&= &\frac{1}{n}\displaystyle\sum_{i=0}^{n-1} d^2(x_{k+i},y)-\frac{1}{n}\displaystyle\sum_{i=0}^{n-1} d^2(x_{k+i},\sigma_n^k),
\end{eqnarray}
which is the intended result. In particular, we have
\begin{equation*}
d^2\big( \sigma_n^k ,y\big)\leqslant \frac{1}{n}\displaystyle\sum_{i=0}^{n-1} d^2(x_{k+i},y).
\end{equation*}
\eqref{l15i2}. Triangle inequality yields:
\begin{equation*}
d^2\big( \sigma_n^k ,y\big)+d^2\big( y,x_{k+i}\big)-2d\big( \sigma_n^k ,y\big)d\big( y,x_{k+i}\big)\leqslant d^2\big( \sigma_n^k ,x_{k+i}\big),
\end{equation*}
hence,
\begin{equation*}
d^2\big( y,x_{k+i}\big)\leqslant d^2\big( \sigma_n^k ,x_{k+i}\big)+2d\big( \sigma_n^k ,y\big)d\big( y,x_{k+i}\big)-d^2\big( \sigma_n^k ,y\big).
\end{equation*}
So summing up over $i$ from $0$ to $n-1$ and multiplying by $\frac{1}{n}$ imply:
\begin{equation}\label{convkarcher1}
\frac{1}{n}\displaystyle\sum_{i=0}^{n-1}d^2\big( y,x_{k+i}\big)\leqslant \frac{1}{n}\displaystyle\sum_{i=0}^{n-1}d^2\big( \sigma_n^k ,x_{k+i}\big)+2d\big( \sigma_n^k ,y\big)\frac{1}{n}\displaystyle\sum_{i=0}^{n-1}d\big( y,x_{k+i}\big)-d^2\big( \sigma_n^k ,y\big).
\end{equation}
On the other hand, by \eqref{ee1} we have:
\begin{equation}\label{convkarcher2}
\frac{1}{n}\displaystyle\sum_{i=0}^{n-1} d^2(x_{k+i},\sigma_n^k)\leqslant \frac{1}{n}\displaystyle\sum_{i=0}^{n-1} d^2(x_{k+i},y)-d^2(\sigma_n^k ,y).
\end{equation}
\eqref{convkarcher1} and \eqref{convkarcher2} show that
\begin{equation*}
d\big( \sigma_n^k ,y\big)\leqslant \frac{1}{n}\displaystyle\sum_{i=0}^{n-1} d(x_{k+i},y).
\end{equation*}
\end{proof}
\begin{lemma}\label{l13}
Let $\mathscr H$ be a Hadamard space and $T:{\mathscr H}\longrightarrow {\mathscr H}$ be a nonexpansive mapping with a nonempty fixed point set $F(T)$. Then for $\{\sigma_n^k (x)\}$ defined in Definition \ref{Karcher mean} and each $k\geqslant 1$ we have:
\begin{enumerate}[(i)]
\item\label{l13i1} The sequence $\{\sigma_n^k(x)\}$ is bounded.
\item\label{l13i2} $d\big(\sigma_n^k(x), T\sigma_n^k(x)\big)\longrightarrow 0$ as $n\rightarrow +\infty$.
\item\label{l13i3} $d\big(\sigma_n(x),\sigma_n^k(x)\big)\longrightarrow 0$ as $n\rightarrow +\infty$.
\item\label{l13i4} $\sigma_n^k(x)\in \overline{co}\{T^mx\}_{m\geqslant 0}$ for all $n\geqslant 1$.
\end{enumerate}
\end{lemma}
\begin{proof}
\eqref{l13i1}.
Since $F(T)\neq \emptyset$, if $p\in F(T)$ by Part \ref{l15i2} of Lemma \ref{l15} and nonexpansiveness of $T$, we see that:
\allowdisplaybreaks\begin{eqnarray}
d\big(\sigma_n^k(x),p\big) &\leqslant & \frac{1}{n}\displaystyle\sum_{i=0}^{n-1} d(T^{k+i}x,p) \nonumber\\
&\leqslant &d(x,p), \nonumber
\end{eqnarray}
thus $\{\sigma_n^k(x)\}$ is bounded.\\
\eqref{l13i2}.
By Part \eqref{l15i1} of Lemma \ref{l15} and  nonexpansiveness of $T$ we have:
\allowdisplaybreaks\begin{eqnarray}
d^2\big(\sigma_n^k(x), T\sigma_n^k(x)\big)&\leqslant & \displaystyle\frac{1}{n}\sum_{i=0}^{n-1}d^2\big(T^{k+i}x, T\sigma_n^k(x)\big)- \displaystyle\frac{1}{n}\sum_{i=0}^{n-1}d^2\big(T^{k+i}x, \sigma_n^k(x)\big)\nonumber \\
&\leqslant &\displaystyle\frac{1}{n}\sum_{i=0}^{n-1}d^2\big(T^{k+i-1}x, \sigma_n^k(x)\big)- \displaystyle\frac{1}{n}\sum_{i=0}^{n-1}d^2\big(T^{k+i}x, \sigma_n^k(x)\big)\nonumber \\
&\leqslant & d^2(T^{k-1}x,\sigma_n^k(x)).\nonumber
\end{eqnarray}
Since the sequences $\{T^nx\}$ and $\{\sigma_n^k(x)\}$ are bounded, the proof is complete.\\
\eqref{l13i3}.
By  Part \eqref{l15i1} of Lemma \ref{l15},  we get:
\allowdisplaybreaks\begin{eqnarray}
d^2\big(\sigma_n(x),\sigma_n^k(x)\big)&\leqslant &\frac{1}{n}\displaystyle\sum_{i=0}^{n-1}d^2\big(T^{k+i}x,\sigma_n(x)\big)- \frac{1}{n}\displaystyle\sum_{i=0}^{n-1}d^2\big(T^{k+i}x,\sigma_n^k(x)\big) \nonumber \\
&\leqslant &\frac{1}{n}\displaystyle\sum_{i=0}^{n-1}d^2\big(T^ix,\sigma_n(x)\big)- \frac{1}{n}\displaystyle\sum_{i=0}^{n-1}d^2\big(T^{k+i}x,\sigma_n^k(x)\big) \nonumber \\
&&+ \frac{1}{n}\displaystyle\sum_{i=0}^{k-1} d^2\big(T^{n+i}x,\sigma_n(x)\big) \nonumber \\
&\leqslant &\frac{1}{n}\displaystyle\sum_{i=0}^{n-1}d^2\big(T^ix,\sigma_n^k(x)\big)- \frac{1}{n}\displaystyle\sum_{i=0}^{n-1}d^2\big(T^{k+i}x,\sigma_n^k(x)\big) \nonumber \\
&&+\frac{1}{n}\displaystyle\sum_{i=0}^{k-1} d^2\big(T^{n+i}x,\sigma_n(x)\big) \nonumber \\
&\leqslant & \frac{1}{n}\displaystyle\sum_{i=0}^{n-1}d^2\big(T^{k+i}x,\sigma_n^k(x)\big)- \frac{1}{n}\displaystyle\sum_{i=0}^{n-1}d^2\big(T^{k+i}x,\sigma_n^k(x)\big) \nonumber\\
&& +\frac{1}{n}\displaystyle\sum_{i=0}^{k-1} d^2\big(T^ix,\sigma_n^k(x)\big)+\frac{1}{n}\displaystyle\sum_{i=0}^{k-1} d^2\big(T^{n+i}x,\sigma_n(x)\big). \nonumber
\end{eqnarray}
Now, since $F(T)\neq \emptyset$ and hence, the sequences $\{T^nx\}$, $\{\sigma_n(x)\}$ and $\{\sigma_n^k(x)\}$ are bounded, we obtain $d\big(\sigma_n(x),\sigma_n^k(x)\big)\rightarrow 0$ as $n\to +\infty$.\\
\eqref{l13i4}.
Let $P:\mathscr H\longrightarrow \overline{co}\{T^mx\}$ be the projection map. On the one hand, by the inequality \eqref{projection} for any $i$ and $n$, we have:
\begin{equation*}
d^2\big(T^{i+k}x,\sigma_n^k(x)\big)\geqslant d^2\big(T^{i+k}x,P\sigma_n^k(x)\big)+ d^2\big(P\sigma_n^k(x), \sigma_n^k(x)\big).
\end{equation*}
On the other hand, by the definition of $\sigma_n^k(x)$, we have:
\begin{equation*}
\displaystyle\frac{1}{n}\sum_{i=0}^{n-1} d^2\big(T^{k+i}x, \sigma_n^k(x)\big) \leqslant \displaystyle\frac{1}{n}\sum_{i=0}^{n-1} d^2\big(T^{k+i}x, P\sigma_n^k(x)\big).
\end{equation*}
Two recent inequalities imply $d^2\big(P\sigma_n^k(x), \sigma_n^k(x)\big)=0$, which is the requested result.\\
\end{proof}
{\bf Proof of Theorem \ref{weakergodickarcher}.}
By a well-known fact for a nonexpansive mapping $T$ in a Hadamard space, $F(T)$ is closed and convex. From the definition of metric projection, we have:
\begin{eqnarray}
d(PT^nx,T^nx)&\leqslant &d(PT^{n-1}x,T^nx) \nonumber\\
&=&d(TPT^{n-1}x,T^nx) \nonumber\\
&\leqslant &d(PT^{n-1}x,T^{n-1}x). \nonumber
\end{eqnarray}
This implies that $\{d(PT^n x,T^nx )\}$ is nonincreasing. By Lemma \ref{the2}, $\{PT^nx \}$ converges strongly to an element $p$ of $F(T)$. The sequence $\{\sigma_n(x)\}$ is bounded by Part \eqref{l13i1} of Lemma \ref{l13}, and hence, by Lemma \ref{lconsubseq} there exists a subsequence $\{\sigma_{n_i}(x)\}$ of $\{\sigma_n(x)\}$ such that $\{\sigma_{n_i}(x)\}$ $\triangle -converges$ to $v\in C$. By Lemma \ref{demiclosedness} and Part \eqref{l13i2} of Lemma \ref{l13} we have $v\in F(T)$. In Lemma \ref{the2}, we see that $\lim_{n\rightarrow\infty} PT^nx=p\in F(T)$, if we show that $v=p$, the proof will be complete. Suppose to the contrary, there is a $\delta>0$ such that
\begin{equation*}
d(p,v)=\delta.
\end{equation*}
By \eqref{projection} and the notation of quasi inner product, it is clear that for all $u\in F(T)$,
\begin{equation*}
\langle T^kxPT^kx,uPT^kx \rangle \leqslant 0.
\end{equation*}
Since $\{d(PT^n x,T^nx )\}$ is nonincreasing, by Lemma \ref{l8} and Cauchy-Schwarz like inequality we have:
\allowdisplaybreaks\begin{eqnarray}
\langle T^kxPT^kx , up\rangle &\leqslant & \langle T^kxPT^kx,PT^kxp \rangle \nonumber \\
&\leqslant & d(T^kx,PT^kx)d(PT^kx,p) \nonumber \\
&\leqslant & d(x,Px)d(PT^kx,p), \nonumber
\end{eqnarray}
so by taking $u=v$, we obtain:
\begin{equation*}
d^2(v,PT^kx)+ d^2(T^kx,p)- d^2(T^kx,v)\leqslant d^2(PT^kx,p)+ 2d(x,Px)d(PT^kx,p).
\end{equation*}
Since $PT^nx\longrightarrow p$ we can choose a positive number $k_0$ such that for any $k\geqslant k_0$
\begin{equation}\label{T^kx}
d^2(T^kx,p)- d^2(T^kx,v)\leqslant 0.
\end{equation}
We know that by $(\overline{Q_4})$ condition, the set $F(p,v)$ is convex. On the other hand, by \eqref{T^kx}, $T^kx\in F(p,v)$. Therefore, by Part \eqref{l13i4} of Lemma \ref{l13}, $\{\sigma_n^{k_0}(x)\}\subseteq F(p,v)$. Also by continuity of the metric function, $F(p,v)$ is closed and hence, by Lemma \ref{lconsubseq} it is $\bigtriangleup-$closed. By Part \eqref{l13i3} of Lemma \ref{l13},  $\sigma_{n_i}^{k_0}(x)$ $\triangle -converges$ to $v$. These facts imply that $v\in F(p,v)$ or equivalently $d(v,p)=0$ i.e., $v=p$, which is a contradiction, and this completes the proof.
\begin{example}
Take $\Bbb{R}^2$ with the river metric which is defined for each $(x,y), (x^\prime,y^\prime)\in \Bbb{R}^2$ as follows.
\begin{equation*}
r\big((x,y),(x^\prime,y^\prime)\big)=
\begin{cases}
|y-y^\prime |, \qquad\qquad\quad\ \    \text{if}\  x=x^\prime \\
|y|+|y^\prime |+|x-x^\prime |, \quad \text{otherwise}
\end{cases}
\end{equation*} 
where $|\cdot|$ denotes absolute value norm. It is known that $\Bbb{R}^2$ with the river metric is an $\Bbb{R}-$tree (and hence a Hadamard space with $(\overline{Q_4})$ condition), which is not locally compact. Take the mapping $T:\Bbb{R}^2\longrightarrow \Bbb{R}^2$ defined by 
\begin{equation*}
T(x,y):=\big(f(x),g(y)\big),
\end{equation*}
for each one to one and nonexpansive mapping $f:\Bbb{R}\longrightarrow\Bbb{R}$ with a nonempty fixed point set and each nonexpansive mapping $g:\Bbb{R}\longrightarrow\Bbb{R}$ with $0\in F(g)$. It is easily seen that $T$ is nonexpansive. Now, using Theorem \ref{weakergodickarcher}, the sequence $\{T^n(x,y)\}$ is $\triangle -mean~convergent$ to a fixed point of $T$. But if $T$ has more than one fixed point (it happens if $f$ or $g$ have two or more fixed points), then the weak ergodic convergence is not concluded from Theorem 2.1 of \cite{Liimatainen2012}.
\end{example}

\section{Weak Ergodic Convergence for Continuous Semigroup of Contractions}
In this section we study the results of the previous section for continuous semigroup of contractions with a nonempty fixed point set and prove the $\triangle -mean~convergence$ of the resulting semigroup to a fixed point of the semigroup in a Hadamard space that satisfies $(\overline{Q_4})$ condition. The proofs are similar to the proofs of the discrete version, therefore we will state only the results without their proofs.\\
Let $C$ be a closed convex subset of a Hadamard space $(\mathscr H,d)$. A 1-parameter continuous semigroup $\mathcal{S}=\{S(t) : t\geqslant 0\}$ of contractions on $C$ is a family of self-mappings $S(t):C\to C$ that satisfy the following conditions:

\begin{enumerate}[(i)]
\item $S(0)x=x,~\text{for}~x\in C;$
\item $S(t+s)x=S(t)S(s)x,~\text{for}~x\in C~\text{and}~t,s\geqslant 0;$
\item $\displaystyle\lim_{t\to t_0}S(t)x=S(t_0)x,~\text{for}~x\in C~\text{and}~t,t_0\geqslant 0;$
\item $d\big(S(t)x,S(t)y\big)\leqslant d(x,y),~\text{for}~x,y\in C~\text{and}~t\geqslant 0.$
\end{enumerate}
Let $F(\mathcal{S})$ denote the common fixed points set of the family $\mathcal{S}$, i.e., \linebreak $F(\mathcal{S})=\bigcap_{t\geqslant 0} F\big(S(t)\big)$. Note that since $F\big(S(t)\big)$ is closed and convex for all $t\geqslant 0$ by \cite[Theorem 12]{Kirk2004}, $F(\mathcal{S})$ is a closed and convex set in a Hadamard space.\\
By a curve in $\mathscr H$ we mean a continuous mapping $c:[0,\infty)\longrightarrow\mathscr H$. $\triangle -convergence$ for a curve has been defined in \cite{bacak2014convex} as follow.
$c$ is said to be $\bigtriangleup-$convergent to a point $x\in\mathscr H$ ($c(t)\overset{\bigtriangleup}{\longrightarrow}x$), if $c(t_n)$ $\triangle -converges$ to $x$ for each sequence $\{t_n\}\subset [0,\infty)$ with $t_n\rightarrow\infty$. 
Also, we can consider a curve $c:[0,\infty)\longrightarrow\mathscr H$ as a net $\{c(t)\}$, and use Definition \ref{weakcondef1} 
for the nets \cite{kakavandi2013}.

\noindent It is clear that the continuous semigroup $\mathcal{S}=\{S(t)\}$ of contractions for a given point $x\in\mathscr H$ represents a curve $x(t):=S(t)x$ and we can use the notion of $\triangle -convergence$ of curve for orbit $\{S(t)x\}$ for any $x\in C$.
\begin{definition}[Karcher mean for a curve and continuous semigroup]\label{Karcher mean for curve}
For a curve $c:[0,\infty)\longrightarrow\mathscr H$, $\sigma_T$ and $\sigma_T^s$(Cesaro mean and Vallee-Poussin mean respect to the Karcher mean on curve $c$), are defined respectively as the unique minimizers of the functions 
\begin{equation*}
\mathcal{G}_T(y)=\frac{1}{T}\displaystyle\int_{0}^{T} d^2(c(t),y)dt,
\end{equation*}
and
\begin{equation*}
\mathcal{G}_T^s(y)=\frac{1}{T}\displaystyle\int_{0}^{T} d^2(c(t),y)dt.
\end{equation*}
Also for an orbit $\{S(t)x\}$, $\sigma_T(x)$ and $\sigma_T^s(x)$(Cesaro mean and Vallee-Poussin mean respect to the Karcher mean on the orbit $\mathcal{S}$), are defined respectively as the unique minimizers of the functions 
\begin{equation*}
\mathcal{G}[x]_T(y)=\frac{1}{T}\displaystyle\int_{0}^{T} d^2(S(t)x,y)dt,
\end{equation*}
and
\begin{equation*}
\mathcal{G}[x]_T^s(y)=\frac{1}{T}\displaystyle\int_{0}^{T} d^2(S(s+t)x,y)dt.
\end{equation*}
 Therefore, Cesaro convergence or the mean convergence for a curve and an orbit of continuous semigroup are defined based on the above means i.e., convergence of $\sigma_T$ or $\sigma_T(x)$.
\end{definition}
It is easy to see that Lemmas \ref{the2}, \ref{l15} and \ref{l13} with similar arguments remain true for a continuous semigroup, and we can rewrite them as follows.
\begin{lemma}
Let $\mathscr H$ be a Hadamard space, $\mathcal{S}=\{S(t) : t\geqslant 0\}$ be a continuous semigroup of contractions, that $F(\mathcal{S})$ is nonempty and $P$ be a metric projection of $\mathscr H$ onto $F(\mathcal{S})$. Then for any $x\in \mathscr H$, $\{PS(t) x\}$ converges strongly to an element $p$ of $F(\mathcal{S})$.
\end{lemma}
\begin{lemma}
Let $c:[0,\infty)\longrightarrow\mathscr H$ be a curve in Hadamard space $\mathscr H$. Then for $\sigma_T^s$ defined as the above, for each $y\in\mathscr H$ and $s\geqslant 0$ we have:
\begin{enumerate}[(i)]
\item
$\displaystyle d^2\big( \sigma_T^s ,y\big)\leqslant \frac{1}{T}\displaystyle\int_{0}^{T} d^2(c(s+t),y)dt-\frac{1}{T}\displaystyle\int_{0}^{T} d^2(c(s+t),\sigma_T^s)dt.$
\item $\displaystyle d\big( \sigma_T^s ,y\big)\leqslant \frac{1}{T}\displaystyle\int_{0}^{T} d(c(s+t),y)dt.$
\end{enumerate}
\end{lemma}
\begin{lemma}
Let $\mathscr H$ be a Hadamard space and $\mathcal{S}=\{S(t) : t\geqslant 0\}$ be a continuous semigroup of contractions with a nonempty common fixed point set $F(\mathcal{S})$. Then for $\{\sigma_T^s (x)\}$ defined by Definition \ref{Karcher mean for curve}, and each $s\geqslant 0$, we have:
\begin{enumerate}[(i)]
\item The sequence $\{\sigma_T^s(x)\}$ is bounded.
\item $d\big(\sigma_T^s(x), S(r)\sigma_T^s(x)\big)\longrightarrow 0$ as $T\rightarrow +\infty$ for each $r\geqslant 0$.
\item $d\big(\sigma_T(x),\sigma_T^s(x)\big)\longrightarrow 0$ as $T\rightarrow +\infty$.
\item $\sigma_T^s(x)\in \overline{co}\{S(t)x\}_{t\geqslant 0}$ for all $T\geqslant 0$.
\end{enumerate}
\end{lemma}
Now we can state the weak ergodic theorem for 1-parameter continuous semigroup of contractions.
\begin{theorem}\label{maintheosemigroup}
Let $(\mathscr{H},d)$ be a Hadamard space that satisfies the $(\overline{Q_4})$ condition, $C$ be a nonempty, closed and convex subset of $\mathscr{H}$ and $\mathcal{S}=\{S(t) : t\geqslant 0\}$ be a continuous semigroup of contractions with $F(\mathcal{S})\neq\emptyset$. Then for each $x\in \mathscr{H}$, the orbit $\{S(t)x\}$ is $\triangle -mean~convergent$ to a common fixed point of $\mathcal{S}$.
\end{theorem}

\begin{example}[Semigroups Generated by Monotone Vector Fields ]
Hadamard manifold is a complete, simply connected Riemannian manifold of nonpositive sectional curvature. Hadamard manifolds are examples of Hadamard spaces \cite{bacak2014convex}\, also, Hadamard spaces with constant curvature satisfy the $(\overline{Q_4})$ condition \cite{kakavandi2013}. If $M$ be a Hadamard manifold and  the mapping $A: D(A)\subset M\rightarrow 2^{TM}$ be a monotone multivalued vector field, by \cite[Theorems 5.1 and 5.2]{IwamiyaOkochi2003}, the Cauchy problem
\begin{equation}\label{eveq1}
\begin{cases}
-x^{\prime}(t)\in Ax(t), &\\
 x(0)=x_0, & 
\end{cases}
\end{equation}
has a global solution that by Proposition 4.2 of \cite{ahmadikhatib2018} is unique. Let $S(t)x_0=x(t)$, then $\mathcal{S}=\{S(t) : t\geqslant 0\}$ is a nonexpansive semigroup (see lemma 4.1 of \cite{ahmadikhatib2018}). It is easy to see that the set of singularities of $A$ (i.e. the set $A^{-1}(0)$) is equal to the set of common fixed points of $\mathcal{S}$. For a complete bibliography and more details about the basic concepts of monotone operators in Hilbert spaces, Hadamard manifolds, the exponential map and also monotone vector fields, we refer the reader to \cite{bauschke2017convex}, \cite{sakai1996riemannian}, \cite{nemeth1999} and \cite{lopezmarquez2009}. Now we have the following result as an application of Theorem \ref{maintheosemigroup}.
\begin{itemize}
\item Let $M$ be a Hadamard manifold of an infinite dimensional with constant curvature and suppose $A: D(A)\subset M\rightarrow 2^{TM}$ is a monotone vector field with at least a singularity point. Then every orbit $\{S(t)x\}$ of the semigroup generated by solutions of \eqref{eveq1} is $\triangle -mean~convergent$ to a singularity of $A$.
\end{itemize}
\end{example}

\bibliographystyle{amsplain}

\bibliography{shortbib}

\end{document}